\documentclass[a4paper,12pt]{article}
\usepackage{amssymb}
\usepackage{amsthm}
\usepackage{amsmath}
\usepackage{amsfonts}

\begin{document}

\def\sect{\section}

\newtheorem{thm}{Theorem}[section]
\newtheorem{cor}[thm]{Corollary}
\newtheorem{lem}[thm]{Lemma}
\newtheorem{prop}[thm]{Proposition}
\newtheorem{propconstr}[thm]{Proposition-Construction}

\theoremstyle{definition}
\newtheorem{para}[thm]{}
\newtheorem{ax}[thm]{Axiom}
\newtheorem{conj}[thm]{Conjecture}
\newtheorem{defn}[thm]{Definition}
\newtheorem{notation}[thm]{Notation}
\newtheorem{rem}[thm]{Remark}
\newtheorem{remark}[thm]{Remark}
\newtheorem{question}[thm]{Question}
\newtheorem{example}[thm]{Example}
\newtheorem{problem}[thm]{Problem}
\newtheorem{excercise}[thm]{Exercise}
\newtheorem{ex}[thm]{Exercise}

\def\Bbb{\mathbb}
\def\cal{\mathcal}
\def\mL{{\mathcal L}}
\def\mC{{\mathcal C}}

\overfullrule=0pt

\def\si{\sigma}
\def\prf{\smallskip\noindent{\it        Proof}. }
\def\call{{\cal L}}
\def\nat{{\Bbb  N}}
\def\la{\langle}
\def\ra{\rangle}
\def\inv{^{-1}}
\def\ld{{\rm    ld}}
\def\trdeg{{tr.deg}}
\def\dim{{\rm   dim}}
\def\th{{\rm    Th}}
\def\rest{{\lower       .25     em      \hbox{$\vert$}}}
\def\ch{{\rm    char}}
\def\zee{{\Bbb  Z}}
\def\conc{^\frown}
\def\acl{{\rm acl}}
\def\cls{cl_\si}
\def\cals{{\cal S}}
\def\mult{{\rm  Mult}}
\def\calv{{\cal V}}
\def\aut{{\rm   Aut}}
\def\ffi{{\Bbb  F}}
\def\ffiti{\tilde{\Bbb          F}}
\def\degs{deg_\si}
\def\calx{{\cal X}}
\def\gal{{\cal G}al}
\def\cl{{\rm cl}}
\def\loc{{\rm locus}}
\def\calg{{\cal G}}
\def\calq{{\cal Q}}
\def\calr{{\cal R}}
\def\caly{{\cal Y}}
\def\aff{{\Bbb A}}
\def\cali{{\cal I}}
\def\calu{{\cal U}}
\def\epsilon{\varepsilon} 
\def\Uu{{\cal U}}
\def\rat{{\Bbb Q}}
\def\ga{{\Bbb G}_a}
\def\gm{{\Bbb G}_m}
\def\cee{{\Bbb C}}
\def\ree{{\Bbb R}}
\def\frob{{\rm Frob}}
\def\Frob{{\rm Frob}}
\def\fix{{\rm Fix}}
\def\Uu{{\cal U}}
\def\proj{{\Bbb P}}
\def\sym{{\rm Sym}}
 
\def\dcl{{\rm dcl}}
\def\calm{{\mathcal M}}

\font\helpp=cmsy5
\def\semdp
{\hbox{$\times\kern-.23em\lower-.1em\hbox{\helpp\char'152}$}\,}

\def\dnfo{\,\raise.2em\hbox{$\,\mathrel|\kern-.9em\lower.35em\hbox{$\smile$}
$}}
\def\dnf#1{\lower1em\hbox{$\buildrel\dnfo\over{\scriptstyle #1}$}}
\def\dfo{\;\raise.2em\hbox{$\mathrel|\kern-.9em\lower.35em\hbox{$\smile$}
\kern-.7em\hbox{\char'57}$}\;}
\def\df#1{\lower1em\hbox{$\buildrel\dfo\over{\scriptstyle #1}$}}        
\def\stab{{\rm Stab}}

\def\si{\sigma}
\def\prf{\smallskip\noindent{\it        Proof}. }
\def\call{{\cal L}}
\def\nat{{\Bbb  N}}
\def\la{\langle}
\def\ra{\rangle}
\def\inv{^{-1}}
\def\ld{{\rm    ld}}
\def\trdeg{{tr.deg}}
\def\dim{{\rm   dim}}
\def\th{{\rm    Th}}
\def\rest{{\lower       .25     em      \hbox{$\vert$}}}
\def\ch{{\rm    char}}
\def\zee{{\Bbb  Z}}
\def\conc{^\frown}
\def\cls{cl_\si}
\def\cals{{\cal S}}
\def\mult{{\rm  Mult}}
\def\calv{{\cal V}}
\def\aut{{\rm   Aut}}
\def\ffi{{\Bbb  F}}
\def\ffiti{\tilde{\Bbb          F}}
\def\degs{deg_\si}
\def\calx{{\cal X}}
\def\gal{{\cal G}al}
\def\cl{{\rm cl}}
\def\loc{{\rm locus}}
\def\calg{{\cal G}}
\def\calq{{\cal Q}}
\def\calr{{\cal R}}
\def\caly{{\cal Y}}
\def\aff{{\Bbb A}}
\def\cali{{\cal I}}
\def\calu{{\cal U}}
\def\epsilon{\varepsilon} 
\def\Uu{{\cal U}}
\def\rat{{\Bbb Q}}
\def\ga{{\Bbb G}_a}
\def\gm{{\Bbb G}_m}
\def\cee{{\Bbb C}}
\def\ree{{\Bbb R}}
\def\frob{{\rm Frob}}
\def\Frob{{\rm Frob}}
\def\fix{{\rm Fix}}
\def\Uu{{\cal U}}
\def\proj{{\Bbb P}}
\def\Ker{{\rm Ker}}

\font\helpp=cmsy5
\def\semdp
{\hbox{$\times\kern-.23em\lower-.1em\hbox{\helpp\char'152}$}\,}

\def\dd{{\rm dd}}
\def\idd{{\rm idd}}
\def\ld{{\rm ld}}
\def\ild{{\rm ild}}
\def\rld{{\rm rld}}
\def\calw{{\cal W}}

\def\vlabel{\label}

\title{An invariant for difference field extensions}
\author{Zo\'e Chatzidakis\thanks{partially supported by MRTN-CT-2004-512234 and by ANR-06-BLAN-0183.}{ (CNRS - Universit\'e Paris 7)}
\and
Ehud Hrushovski\thanks{thanks to the Israel Science Foundation (1048/07) for
  support, and to the ANR-06-BLAN-0183.}\  (The Hebrew University)}
\date{}
\maketitle

\section*{Introduction} 
A {\it difference field} is a field with a
distinguished endomorphism $\si$. 
In this short note, we introduce a new invariant for finitely generated
difference field 
extensions of finite transcendence degree, the {\em distant degree}. If
$(K,\si)$ is a difference field, and $a$ a finite tuple in some difference
field extending $K$, and which satisfies $\si(a)\in K(a)^{alg}$ (the
field-theoretic algebraic closure of $K(a)$), we define
$$\dd(a/K)=\lim_{k\rightarrow +\infty}[K(a,\si^k(a)):K(a)]^{1/k}.$$
One shows easily that $\dd(a/K)$ is bounded by a classical invariant of
difference field extensions, 
the {\em limit degree} of $a$
over $K$, and which is defined by  $$\ld(a/K)=\lim_{k\rightarrow
  +\infty}[K(a,\si(a),\ldots,\si^{k+1}(a)):K(a,\si(a),\ldots,
\si^k(a))].$$ 
Our main result, Theorem \ref{sc-thm2}, is that  there is some
$b\in K(a)_\si$ (the difference field generated by $a$ over $K$) such
that $a\in K(b)^{alg}$, and $\dd(b/K)=\ld(b/K)$. In
characteristic $0$, this result is a consequence of a result of
George Willis on scale functions of automorphisms of totally disconnected locally compact groups,
see \cite{W1}, \cite{W2}.

Theorem \ref{sc-thm2} follows immediately from Theorem \ref{sc-thm1},
which asserts that there is $b\in K(a)_\si$ such that $a\in K(b)^{alg}$
and $\si(b)\in K(b,\si^\ell(b))$ for every $\ell>0$. This latter result
is particularly useful for difference fields - it is quite convenient to
find a tuple satisfying $[K(a,\si^\ell(a)):K(a)]=\ld(a/K)^\ell$ for all
$\ell>0$.  We then proceed to
derive other properties of these tuples $b$ satisfying ``ld=dd'', see
Proposition \ref{prop2}. We conclude the study of $\dd$ with Proposition
\ref{prop1}, which among other things shows that $\dd(a,b/K)\geq
\dd(a/K(b)_\si)\dd(b/K)$. Unfortunately, the distant degree is not
multiplicative in towers (see \ref{example}).


The above results continue to hold for the class of perfect fields, in place of the class of fields.  More generally, the statements and proof
go through verbatim for {\em strongly minimal sets},
cf. e.g. \cite{pillay} for a definition.  Fields should be replaced by definably closed substructures $K$ of a model 
$M$ of the given strongly minimal theory.  We then obtain an invariant of automorphisms of such substructures.

The results for strongly minimal sets admit a purely group theoretic
presentation.  Namely let  $G$ be a group, $\si$  an automorphism of
$G$, and $H$  a
subgroup  of $G$ such that $H^\si\cap H$ has finite index in $H$ and  in $H^\si$.  Then one can define the distant degree in terms of $(G,H,\si)$ alone.
When  $\calu$ is a strongly minimal structure with an automorphism $\si$, $K$ a substructure, $a \in \calu \setminus K$, setting  $G=\aut(\calu/K)$,  
$H=\aut(\calu/K(a))$, and $H^\si=\aut(\calu/K(\si(a)))$, we recover the
previous definitions.   See the earlier ArXiv version of the paper for
details, {\tt http://arxiv.org/abs/0902.0844v2}. 

After formulating the results group-theoretically, we found earlier
results of Willis extending most of ours in this
context \footnote{thanks to Dugald Macpherson for drawing Willis's
  results to our attention}.  Willis starts out from a 
totally disconnected locally compact
group, rather than an abstract group $G$ with a subgroup $H$ as above;  one can however complete the abstract group $G$ above with respect to the topology generated by the finite index subgroups of $H$; so
again the two settings are equivalent. It follows that our invariant
$\dd(a/K)$ coincides with the {\em scale} of $\si$ in the sense of
Willis. This  yields two new ways of computing the scale function: the
definition of dd, and 
 Lemma \ref{sc-lem1}(3).

Willis' results allowed us to strengthen our original results. A key
observation towards Theorem \ref{sc-thm1}
comes from a result
hidden in  Lemma 3(a) of \cite{W1}. Further help  comes from the definition
of Willis' group $\call$, but the other ingredients in our  proof are 
different.

We conclude the paper with a discussion of the three settings.  In \ref{Wil1}
- \ref{Wil3} we compare our results in the field setting with Willis' in the group
setting; naturally they bring in intuitions from different
directions. We then show the
equivalence of the setting of strongly minimal structures with the one
of 
 totally disconnected locally compact groups, see \ref{sm}.


At the end of chapter 1, we also refine the main  results for definable groups.  By a {\em difference subgroup} we mean here
 a subgroup of an algebraic group defined by difference equations; by a
{\em morphism}, we mean a group homomorphism
 given locally in the $\si$-topology by difference-rational functions. We show
 in Proposition \ref{prop3} that  if $H$ is a difference subgroup,
 has finite order and is connected for the $\si$-topology, then  there is a 
morphism $f:H\to H'$ with finite central kernel, such that if $b$ is a generic
of the difference subgroup $H'$, then $\ld(b/K)=\dd(a/K)$.

%


\section{The results}
\para{\bf Setting, notation and convention}. A {\it difference field} is a field with a
distinguished endomorphism $\si$. If $\si$ is onto, it is called an {\it
inversive} difference field. Every difference field $(K,\si)$ has an
{\it inversive closure}, denoted $K^{inv}$, which is characterised by
admitting a unique $K$-embedding
into any inversive difference field containing $K$ ([Co], 2.5.II). We will work in
some large inversive difference field $(\calu,\si)$.

If $a$ is a tuple in $\calu$, then $K(a)_\si$ denotes the difference
field generated by $a$ over $K$, i.e., $K(a)_\si=K(\si^i(a)\mid
i\in\nat)$. 
If $E$ is a field, then $E^{alg}$ denotes the (field-theoretic)
algebraic closure of $E$, $E^s$ its separable closure, and $E^{perf}$
its perfect hull. If $a$ is a tuple in $E^{alg}$, then
$\mu(a/E)$ denotes $[E(a):E]$. 

We will say that a sequence $(a_n)_{n\in\nat}$ is {\em increasing} if $a_n\leq
a_{n+1}$ for any $n\in\nat$. Similarly for decreasing.

\para{\bf Definitions}. Let $K$ be a difference subfield of $\calu$, $a$
be a finite tuple in $\calu$, and assume that 
$\si(a)\in K(a)^{alg}$. 
\begin{enumerate}
\item{The {\it limit degree of $a$ over $K$} (or of $K(a)_\si$ over
$K$) is $$\ld(a/K)=\lim_{k\rightarrow
\infty}\mu(\si^{k+1}(a)/K(a,\si(a),\ldots,\si^k(a))),$$
and the {\it inverse limit degree of $a$ over $K$} is
$$\ild(a/K)=\lim_{k\rightarrow
\infty}\mu(\si^{-(k+1)}(a)/K^{inv}(a,\si\inv(a),\ldots,\si^{-k}(a))).$$
}
\item{We
define the {\it distant degree} and {\it inverse distant degree of $a$
  over $K$}  by $$\dd(a/K)=\lim_{k\rightarrow 
\infty}\mu(\si^{k}(a)/K(a))^{1/k},\quad \idd(a/K)=\lim_{k\rightarrow
\infty}\mu(\si^{-k}(a))/K^{inv}(a))^{1/k}.$$}
\end{enumerate}

\para\vlabel{propld}{\bf Properties of the limit degree}.
The limit and inverse limit degrees  are invariants of the
extension $K(a)_\si/K$, they are multiplicative in towers, and
$\ld(a/K)=\ld(a/K^{inv})$, see [Co], section 5.16. If  
$\mu(\si(a)/K(a))=\ld(a/K)$, then for every $i\in\nat$, the fields
$K(\si^j(a)\mid j\geq i)$ and $K^{inv}(\si^j(a)\mid j\leq i)$ are linearly
disjoint over $K(\si^i(a))$. Indeed, the numbers
$\mu(\si^k(a)/K(a,\ldots,\si^{k-1}(a)))$ form a decreasing sequence, and
$\ld(a/K)$ is the value at which it stabilises. Thus, when
$\mu(\si(a)/K(a))=\ld(a/K)$, 
$\ld(\si^i(a)/K)=\mu(\si^{i+1}(a)/K(\si^i(a)))$ for every $i\geq 0$. From
$$\mu(\si^{i+1}(a)/K(\si^i(a)))=\ld(\si^i(a)/K)=\ld(\si^i(a)/K^{inv})=\mu(\si^{i+1}(a)/K^{inv}(\si^j(a)\mid
j\leq i)),$$
we obtain that $K^{inv}(\si^j(a)\mid
j\leq i))$ and $K(\si^i(a),\si^{i+1}(a))$ are linearly disjoint over
$K(\si^i(a))$. An easy induction argument gives the result. 
In this case one also has 
$\ild(a/K)=\mu(a/K^{inv}(\si(a)))$. Furthermore, if $i<j<k$,
then $$\mu(\si^j(a)/K^{inv}(\si^i(a),\si^k(a)))=
\mu(\si^j(a)/K^{inv}(\si^\ell(a), \ell\in (-\infty,i]\cup [k,+\infty))).\eqno{(\#)}$$

\para\vlabel{lem4}{\bf Lemma}. Let
$a$ and $b$ be tuples in $\calu$ such that $b,\si(a)\in K(a)^{alg}$,
$\si(b)\in K(b)^{alg}$.
\begin{itemize}
\item [(1)]
There is a constant $D$ such that for all
$k\in\nat$, $\mu(\si^k(a),\si^k(b)/K(a,b))\leq
D\mu(\si^k(a))/K(a))$. Hence $\dd(b/K)\leq \dd(a/K)$. 

\item[(2)] There is a constant $D'$ such that for every $k>0$,
$\mu(\si^k(a)/K(a))\leq D'\mu(\si^k(a)/K^{inv}(a))$.
\item [(3)] $\ld(a,b/K)\ild(a/K)=\ild(a,b/K)\ld(a/K)$.

\end{itemize}

\prf (1) One verifies easily that 

\begin{eqnarray*}
\mu(\si^k(a),\si^k(b)/K(a,b))&\leq& 
\mu(\si^k(b)/K(\si^k(a)))\mu(\si^k(a)/K(a))\cr
&&\qquad\leq \mu(b/K(a))\mu(\si^k(a)/K(a)).
\end{eqnarray*}
Take
$D=\mu(b/K(a))$. 

(2) Let $n$ be such that
$\mu(\si^{n+1}(a)/K(a,\ldots,\si^n(a)))=\ld(a/K)$, and let
$E=K(a,\ldots,\si^n(a))$. Then, for $m\geq n$,
we have 
$\mu(\si^m(a)/K(a))\leq \mu(\si^m(a)/E)\mu(E/K(a))$, and
$\mu(\si^m(a)/E)=\mu(\si^m(a)/K^{inv}(E))\leq
\mu(\si^m(a)/K^{inv}(a))$. Take $D'=\mu(E/K(a))$.

(3) Using the multiplicativity in towers of the limit degrees and inverse limit
degrees the desired equality becomes
$$\ld(a/K)\ld(b/K(a)_\si)\ild(a/K)=\ild(a/K)\ild(b/K(a)_\si)\ld(a/K).$$
Therefore (and using (2)), it suffices to show that if $L=K(a)_\si^{inv}$
then 
$\ld(b/L)=\ild(b/L)$. We have
$\mu(b/L)=\mu(\si(b)/L)$, so that 
$$\mu(\si(b)/L(b))={\mu(b,\si(b)/L)\over \mu(b/L)}=\mu(b/L(\si(b))).$$
If $\ld(b/L)=\mu(\si(b)/L(b))$, this gives the result. Else, it
suffices to replace $b$ by $(b,\si(b),\ldots,\si^n(b))$ for some $n$.

\para{\bf Setting}. The previous lemma has three immediate consequences:
if $K(a)^{alg}_\si=K(b)_\si^{alg}$, then $\dd(a/K)=\dd(b/K)$ (item 1);
$\dd(a/K)=\dd(a/K^{inv})$ (item 2); and
$\dd(a/K)\ild(a/K)=\idd(a/K)\ld(a/K)$ (item 3). This 
reduces the study of $\dd$ to the following setting: we work inside a
large algebraically closed 
difference field $\calu$, over a difference field $K=\si(K)$, and $a$ is a tuple
such that $\si(a)\in K(a)^{alg}$ and $\mu(\si(a)/K(a))=\ld(a/K)$.

\para\vlabel{sc-lem1}{\bf Lemma}. \begin{enumerate}
\item[(1)]The sequence $\mu(\si(a)/K(a,\si^\ell(a)))$, $\ell\in\nat^{>0}$, is
  an increasing sequence.
\item[(2)]Let $m=\sup\{\mu(\si(a)/K(a,\si^\ell(a))),\ \ell\in\nat^{>0}\}$,
  let $\ell_0$ be the smallest $\ell$ at which this value is
  attained, and let $C=\mu(\si(a),\ldots,\si^{\ell_0-1}(a)/a,\si^{\ell_0}(a))$. If $\ell,j\geq \ell_0$, then
$$\mu(a/K(\si^{-j}(a),
\si^\ell(a)))={m^{\ell_0}\over C}.$$
\item[(3)] With $m$ as in (2), $$\dd(a/K)={\ld(a/K)\over m}.$$
\end{enumerate}

\prf We will omit $K$ from the notation, i.e., $\mu(a/b)$ denotes
$\mu(a/K(b))$. We will use equation $(\#)$ of \ref{propld} repeatedly. 

(1) One
has 
$$\mu(\si(a)/a,\si^{\ell}(a))=\mu(\si(a)/a,\si^\ell(a),\si^{\ell+1}(a))\leq
\mu(\si(a)/a,\si^{\ell+1}(a)).$$
(The first equality is an example of the use of  \ref{propld} $(\#)$).

(2) If $\ell\geq \ell_0$, then\\[0.1in]
$\mu(\si(a),\ldots,\si^{\ell-1}(a)/a,\si^\ell(a))=\prod_{j=1}^{\ell-1}\mu(\si^i(a)/K(\si^{i-1}(a),\si^\ell(a)))$\\
\hspace*{3cm}$=\prod_{j=1}^{\ell-\ell_0}\mu(\si^j(a)/\si^{j-1}(a),\si^{\ell}(a))
\mu(\si^{\ell-\ell_0+1}(a),\ldots,\si^{\ell-1}(a)/\si^{\ell-\ell_0}(a), 
  \si^\ell(a))$\\
\hspace*{3cm}$=m^{\ell-\ell_0}C.$

If $j\geq\ell_0$, applying $\si^{-j}$ to the above equation with
$\ell=j$ gives
$\mu(\si^{-j+1}(a),\ldots,\si\inv(a)/\si^{-j}(a),a)=m^{j-\ell_0}C$. 

On the other hand, 
\begin{eqnarray*}\mu(\si^{-j+1}(a),\ldots,\si^{\ell-1}(a)/\si^{-j}(a),\si^\ell(a))&=&\mu(a/\si^{-j}(a),\si^\ell(a))\mu(\si^{-j+1}(a),\ldots,\si^{-1}(a)/\si^{-j}(a),a)\\ 
&&\qquad\qquad\qquad\mu(\si(a),\ldots,\si^{\ell-1}(a)/a,\si^\ell(a))\\
&=&\mu(a/\si^{-j}(a),\si^\ell(a))Cm^{j-\ell_0}Cm^{\ell-\ell_0},
\end{eqnarray*}
which implies that
$$\mu(a/\si^{-j}(a),\si^\ell(a))={Cm^{j+\ell-\ell_0}\over
  C^2m^{j+\ell-2\ell_0}}={m^{\ell_0}\over C}.$$

(3) We computed in the proof of (2) that for $\ell\geq \ell_0$,
$\mu(\si(a),\ldots,\si^{\ell-1}(a)/a,\si^\ell(a))=Cm^{\ell-\ell_0}$. Hence
$$\mu(\si^\ell(a)/a)={\mu(\si(a),\ldots,\si^\ell(a)/a)\over
  \mu(\si(a),\ldots,\si^{\ell-1}(a)/a,\si^\ell(a))}={\ld(a/K)^\ell\over
  Cm^{\ell-\ell_0}}=\left({\ld(a/K)\over m}\right )^\ell {m^{\ell_0}\over
  C}.$$

\para{\bf Definition}. \vlabel{minpol} Let $a=(a_1,\ldots,a_n)$ be algebraic over the
field $L$. We define the {\em tuple of minimal monic polynomials of $a$
  over $L$} as follows: $p=(p_1,\ldots,p_n)$, with $p_i\in L[X_1,\ldots,X_i]$,
$i=1,\ldots,n$, are 
such that $p_1(X_1)$ is the minimal monic polynomial of $a_1$ over $L$, and
for $1<i\leq n$, $p_i(a_1,\ldots,a_{i-1},X_i)$ is the minimal monic polynomial
of $a_i$ over $L(a_1,\ldots,a_{i-1})=L[a_1,\ldots,a_{i-1}]$. Then
$\mu(a/L)=\prod_{i}\deg_{X_i}p_i$. 

Let $L_0$ be a subfield of $L$, and assume that
$\mu(a/L_0)=\mu(a/L)$. Then the tuple $p$ has its coefficients in
$L_0$. This follows from the fact that for any subfield $L_0$ of $L$,
one always has $\mu(a_i/L(a_1,\ldots,a_{i-1}))\leq \mu
(a_i/L_0(a_1,\ldots,a_{i-1}))$ for $i=1,\ldots,n$, so that our
assumption on the degree of the extension forces equality
everywhere. 

\para\vlabel{sc-thm1}{\bf Theorem}. Let $K=\si(K)$, and $a$ a tuple such
that $\si(a)\in 
K(a)^{alg}$. Then there is $c\in K(a)_\si$ such that $a\in K(c)^{alg}$,
and for every $\ell>i>0$, $\si^i(c)\in K(c,\si^\ell(c))$. 

\prf We may assume that $\mu(\si(a)/K(a))=\ld(a/K)$. We let $\ell_0$,
$m$ and $C$ be defined as in Lemma \ref{sc-lem1}, and let $c$ be the tuple
of coefficients of the tuple of minimal monic polynomials of $a$ over
$K(\si^{-\ell_0}(a),\si^{\ell_0}(a))$. 

Since
$\mu(\si(a)/K(a))=\ld(a/K)$, we have $\mu(a/K(\si^i(a)\mid |i|\geq
\ell_0))=\mu(a/K(\si^{-\ell_0}(a),\si^{\ell_0}(a)))$. Hence, using
Lemma \ref{sc-lem1}, if $j,\ell\geq\ell_0$, then $c$ belongs to
$K(\si^{-j}(a),\si^\ell(a))$. Let $$F=\bigcap_{\ell-n\geq 2\ell_0}K(\si^i(a)\mid i\in
(-\infty,n]\cup [\ell,+\infty)).$$   Then $c\in F$ and $\si(F)=F$. We have $\mu(a/F)=\mu(a/K(c)):=N$.
Let $\ell\geq \ell_0$. Then 
$\mu(\si^{-\ell}(a)/F(a, \si^\ell(a)))=N$ because $F(a,
    \si^\ell(a))\subseteq K(\si^i(a)\mid i\in
(-\infty,-\ell-\ell_0]\cup [-\ell+\ell_0,+\infty))$ and
$\si^{-\ell}(c)\in F$; and 
$\mu(\si^\ell(a)/F( a))=N$ because $F(a)\subseteq K(\si^i(a)\mid i\in
(-\infty,\ell-\ell_0]\cup [\ell+\ell_0,+\infty))$ and $\si^{\ell}(c)\in F$.
 This implies that
$$[K(\si^{-\ell}(a),\si^\ell(a),
\si^{-\ell}(c),c,\si^\ell(c)):K(\si^{-\ell}(c),c,\si^\ell(c))]=N^2,$$
    and therefore that $$c\in K(\si^{-\ell}(c),\si^\ell(c)).$$ 
The first implication is clear; for the second, we know that $c$
belongs to $K(\si^{-\ell}(a),\si^\ell(a))$, so if $c\notin
K(\si^{-\ell}(c),\si^\ell(c))$, we would have
$\mu(\si^{-\ell}(a),\si^\ell(a)/K(\si^{-\ell}(c),c,\si^\ell(c)))<N^2$.

Assume that $\si(c)\notin K(c,\si^\ell(c))$ for some $\ell>0$, and let
$n$ be the maximum  value of $\mu(\si(c)/K(c,\si^{\ell}(c)))$, attained at
$\ell_2$ but not before. As we saw in Lemma \ref{sc-lem1},
if $\ell\geq\ell_2$ and  
$C':=\mu(\si(c),\ldots,\si^{\ell_2-1}(c)/K(c,\si^{\ell_2}(c)))$, then
$\mu(c/K(\si^{-\ell}(c),\si^\ell(c)))=n^{\ell_2}/C'$,  
i.e.,
$n^{\ell_2}=C'$
(since for $\ell\gg 0$, $c\in K(\si^{-\ell}(c),\si^\ell(c))$). But
by definition of $\ell_2$, if $j<\ell_2$, then
$\mu(\si(c)/c,\si^j(c))<n$. Hence 

$$C'=\prod_{i=1}^{\ell_2-1}\mu(\si^i(c)/\si^{i-1}(c),\si^{\ell_2}(c))=n^{\ell_2},$$ 
which implies $n=1$, since the second term is $\leq
n^{\ell_2-1}$. I.e., $\si(c)\in K(c,\si^\ell(c))$ for all $\ell>
0$. An easy induction then gives that $\si^i(c)\in K(c,\si^\ell(c))$ if
$0<i<\ell$. The proof gives that $c\in K(a)_\si^{inv}$;  if $m$ is such that
$\si^m(c)\in K(a)_\si$, then $\si^m(c)$ is our desired element. 

\para \vlabel{sc-thm2} We will now derive some consequences of Theorem
\ref{sc-thm1}. First note a 
very easy corollary:

\smallskip\noindent
{\bf Theorem}. Let $K=\si(K)$, $a$ such that $\si(a)\in
K(a)^{alg}$, and let $c$ be given by Theorem \ref{sc-thm1}. Then
$\dd(a/K)=\ld(c/K)$.

\prf By Lemma \ref{lem4}, $\dd(a/K)=\dd(c/K)$. On the other hand, since
$\si(c)\in K(c,\si^\ell(c))$ for every $\ell>0$, we have $\mu(\si^\ell(c)/K(c))=
\ld(c/K)^\ell$. \\

We now proceed to list properties of elements satisfying $\ld=\dd$.

\para\vlabel{prop2}{\bf Proposition}. Let $K=\si(K)$, $a$ a tuple such that $\si(a)\in
K(a)^{alg}$, and $c\in K(a)_\si$ given by Theorem \ref{sc-thm1}.
\begin{itemize}
\item[(1)] The following conditions are equivalent, for a tuple $d$
  which is equi-algebraic with $a$ over $K=\si(K)$:
\begin{itemize}
\item[(i)] $\ld(d/K)=\dd(a/K)$ ($=\dd(d/K)$).

\item[(ii)]$\ld(d/K)=\inf \{\ld(e/K)\mid K(e)^{alg}=K(a)^{alg}\}$.
\end{itemize}
If in addition $\mu(\si(d)/K(d))=\ld(d/K)$, then each of the above
conditions is equivalent to each of the following:
\begin{itemize}
\item[(iii)] For every $\ell>0$, $\si(d)\in
  K(d,\si^\ell(d))$. 
\item[(iv)] For every $\ell>0$, $d\in K(\si^{-\ell}(d),\si^\ell(d))$.
\end{itemize}
Furthermore, any of the above conditions is equivalent to the analogous
one for $\si\inv$.
\item[(2)] 
  Let $b$ be the set of conjugates of $a$ over $K(c)_\si$, and let $d$ 
 be a code for the set $b$ (i.e., $K(c)_\si(d)$ is the subfield of
 $K(c)_\si(b)$ fixed under $\aut(K(c)_\si(b)/K(c)_\si)$). Then for some
 $n$, $\ld(d^{p^n}/K)=\dd(d^{p^n}/K)=\dd(a/K)$. If $K$ is perfect, then
  $\ld(d/K)=\dd(a/K)$, and $a\in K(d)^s$. 
\item[(3)] The number
  $\dd_{\si^n}(a/K)$ computed in the $\si^n$-difference field $\calu$,
  equals the $n$-th power of $\dd_\si(a/K)$. 
\item[(4)] $\dd(a/K)=1$ if and only if $\{\mu(\si^\ell(a)/K(a))\mid
  \ell\in\nat\}$ is bounded. In that case, $\si(c)\in K(c)$.
\item[(5)] $\dd(a/K)$ divides $\ld(a/K)$. 
\item[(6)] Assume that $\ld(d/K)=\dd(d/K)$. Then also
  $\ld(c,d/K)=\dd(c,d/K)$. 
\item[(7)] Assume that $d$ is equi-algebraic with $a$ over $K$, and that
  for some $\ell_1$, $d\in \bigcap_{\ell\geq
    \ell_1}K(\si^{-\ell}(a),\si^\ell(a))$. Then $\dd(d/K)=\ld(d/K)$. 
\end{itemize}

\prf (1) The limit degree satisfies
$\ld(a/K)=\ld(a,\si(a),\ldots,\si^n(a)/K)$ for every $n$, and we may
therefore assume that $\ld(a/K)=\mu(\si(a)/K(a))$ since this change will
not affect the first two conditions. We will show the equivalence of (i)
-- (iv).

We know by Lemma \ref{lem4} that $\dd(a/K)=\dd(d/K)$. Assume
that (iii) does not hold. Then for some $\ell>0$, we have $\si(d)\notin
K(d,\si^\ell(d))$; by Lemma \ref{sc-lem1} (1) and (3), 
we have 
 $\dd(d/K)<\ld(d/K)$, whence $\dd(a/K)<\ld(d/K)$. Thus (i) implies
(iii). Clearly (iii) implies (i). 

Similarly,
$\dd(e/K)\leq \ld(e/K)<\dd(a/K)$ is impossible unless $K(e)^{alg}$ is strictly
contained in $K(a)^{alg}$, and this proves the equivalence of (i) and
(ii).

(iii) implies (iv) is an easy induction, and (iv) implies (iii) is proved
in the last part of the proof of Theorem \ref{sc-thm1}.

Finally, for the last assertion it suffices to show that one of the
above conditions is equivalent to 
its analogue for $\si\inv$. 
We know that the quotient ${\ld(a/K)\over \ild(a/K)}$ is an invariant
  of the extension $K(a)^{alg}/K$, by Lemma \ref{lem4}(3). Hence, (ii) for $\si$ implies (ii)
  for $\si\inv$. 

(2) By
definition of $d$ the extension $K(b)/K(d)$ is separable, and
the extension $K(c)_\si(d)/K(c)_\si$ is purely inseparable. This implies
that $K(c,d)_\si/K(c)_\si$ is purely inseparable, and  $a\in K(d)^s$.  If
$n$ is such that $d^{p^n}\in 
K(c)_\si$, then 
$\ld(d^{p^n}/K)$ divides $\ld(c/K)$, and by minimality of the
latter, must be equal to it. Hence $\ld(d^{p^n}/K)=\dd(a/K)$. 

If $K$ is perfect, then
$\ld(d^{p^n}/K)=\ld(d/K^{p^{-n}})=\ld(d/K)=\dd(a/K)$. 
 
(3) Clear from the definition of $\dd$.

(4) Clear by Lemma  \ref{lem4}(1) and Theorem \ref{sc-thm2}. 

(5) As $c\in K(a)_\si$, $\dd(a/K)=\ld(c/K)$ divides $\ld(a/K)$. 

(6) By (1), we have $\si(c)\in K(c,\si^\ell(c))$ and $\si(d)\in
K(d,\si^\ell(d))$ for every $\ell>0$. Hence, $\si(c,d)\in
K(c,d,\si^\ell(c,d))$ for every $\ell>0$, which by (1) implies that
$\ld(c,d/K)=\dd(c,d/K)$. 

(7) We use the notation of Theorem \ref{sc-thm1}. Without loss of
generality, we may assume that $\ell_1\geq \ell_0$. Let $e$ be a tuple
such that $K(e)=\bigcap_{\ell\geq
  \ell_1}K(\si^{-\ell}(a),\si^\ell(a))$. Then $c\in K(e)$ (since
$\ell_1\geq \ell_0$), and therefore is equi-algebraic with $e$ over
$K$. As $d\in K(e)$, it suffices to show that $\ld(e/K)=\dd(e/K)$, since
$\ld(d/K)\leq \ld(e/K)$, and by (1). 

Let $F_0$ be the inversive difference field generated by $K(e)$. Then
$F_0\subseteq F$, and $c\in F_0$. These imply that
$\mu(\si^{-\ell}(a)/F_0(a,\si^\ell(a)))=N=\mu(\si^\ell(a)/F_0(a))$. Reasoning
as in the proof of Theorem \ref{sc-thm1} one gets $e\in
K(\si^{-\ell}(e),\si^\ell(e))$. Now use (1) to conclude. \\

We now investigate the behaviour of $\dd$ in towers of
extensions. Unfortunately, it is not multiplicative, as we will see in
\ref{example}.

\para\vlabel{prop1}{\bf Proposition}. Let $K\subset \calu$ be a
difference field, $a$ and $b$ two tuples in $\calu$ such that $\si(a)\in
K(a)^{alg}$, $\si(b)\in K(b)^{alg}$.
\begin{enumerate}
%
\item{ $\dd(a,b/K)\geq \dd(a/K(b)_\si)\dd(b/K)$.}
\item{If $b\in K(a)^{alg}$, then $\dd(b/K)\leq \dd(a/K)$.}
\end{enumerate}

\prf 
(1) By 
Lemma \ref{lem4}(2) we may assume that $K$ is
inversive. Let $d$ be a finite tuple of $K(b)^{alg}$ such that
$K(a,b,d)$ is a regular extension of $K(b,d)$. If $C=[K(b,d):K(b)]$,
then 
 for any $\ell>0$,
$$\mu(\si^\ell(b)/K(b))\leq C\mu(\si^\ell(b)/K(a,b)).$$
Thus
    \begin{eqnarray*}\mu(\si^\ell(a),\si^\ell(b)/K(a,b))&=&
    \mu(\si^\ell(a)/K(a,b,\si^\ell(b)))\mu(\si^\ell(b)/K(a,b))\cr
&\geq& \mu(\si^\ell(a)/K(b)_\si(a))C\inv\mu(\si^\ell(b)/K(b)).
\end{eqnarray*}
This gives the result.

(2) Follows immediately   Lemma \ref{lem4}.

\para\vlabel{example}{\bf An example}. Unfortunately, Proposition
\ref{prop1}(1) is the best we can hope for, 
the invariant $\dd$ is not multiplicative in towers. Here is an example.

\smallskip
Let $a$ be a generic solution of $\si(a^2)=a^2+1$ over  an algebraically
closed inversive difference field $K$ of characteristic $0$, and
$b$ a solution of $\si(b)=b+a$. Then $\dd(a/K)=\dd(a^2/K)=\ld(a^2/K)=1$, $\ld(a/K)=2$,
and $\ld(b/K(a)_\si)=1=\dd(b/K(a)_\si)$, so that $\ld(a,b/K)=2$. If $\ell>0$, then
$\si^\ell(b)-b= a+\sqrt{a^2+1}+\cdots +\sqrt{a^2+\ell-1}$,
$K(a^2,\si(b)-b)=K(a,\sqrt{a^2+1},\ldots,\sqrt{a^2+\ell-1})$ is an
extension of degree $2^\ell$ of $K(a^2)$. 

Thus, if $\ell>1$, then
$\si(a), \si(b)\in K(a,b, \si^\ell(a),\si^\ell(b))$, so that
$\dd(a,b/K)=\ld(a,b/K)=2$, but $\dd(a/K)\dd(b/K(a)_\si)=1$. 

\para{\bf Remark}. Note that the example shows that the failure of
multiplicativity in towers is fundamental: taking
$L=K(a)^{alg}$ and $M=K(a,b)^{alg}$, we obtain a tower $K\subset
L\subset M$ of algebraically closed inversive difference fields with
$$\dd(M/K)=2\neq \dd(L/K)\dd(M/L).$$

\para\vlabel{lem5}{\bf The case of difference subgroups of algebraic groups}. In case
our tuple $a$ is the 
generic of some difference subgroup, we will show that the tuple $c$ can be
chosen to be 
the generic of a difference  subgroup, with the map $a\mapsto c$ a 
morphism. We first need a lemma:

\smallskip\noindent
{\bf Lemma}. Let $K$ be a field,
$G_1,G_2, U$  
(connected) algebraic groups defined over $K$ with 
$U\subset G_1\times G_2$, and $\pi_i:G_1\times G_2\to G_i$ the natural
projections. Assume that $\pi_i(U)=G_i$ for $i=1,2$. If $S_1=\pi_1(U\cap
(G_1\times 1))$, then 
$S_1$ is a normal subgroup of $G_1$. Moreover, if the restriction of
$\pi_2$ to $U$ is finite, if $g=(a,b)$ is a generic of
$U$ over $K$, then the field conjugates of $a$ over
$K(b)$ are the elements of $a+S_1$, and $S_1$ is central. 

\prf $S_1\times 1=\Ker(\pi_2)$ is
normal in $U$, and because $\pi_1(U)=G_1$, $S_1$ is normal in $G_1$. The
finiteness of $\pi_2\rest_U$ implies that $S_1$ is finite, and therefore
central since $G_1$ is connected. If $a'\in
a+S_1$, then $(a',b)$ is also a generic of $U$ over $K^{alg}$, and
therefore the fields $K(a,b)$ and $K(a',b)$ are
$K(b)$-isomorphic.

\para \vlabel{prop3}
{\bf Proposition}. Assume that $K$ is a difference
field, let $H$ be a difference subgroup of some algebraic group $G$, both
defined over $K$, and
assume that if $a$ is a generic of $H$ (for the $\si$-topology), then
$K(a)_\si$ is a regular extension of $K$ of finite transcendence degree
over $K$. Then there are a difference subgroup $H'$, 
a morphism $f:H\to H'$ with finite central kernel, defined by a
tuple of  difference rational 
functions, and such that  $\ld(f(a)/K)=\dd(a/K)$. 

\prf Let $a$  be a generic of $H$ over $K$. Choose $n$ such that
$\si^{n+1}(a)\in K(a,\ldots,\si^n(a))^{alg}$ and 
$\mu(\si^{n+1}(a)/K(a,\ldots,\si^n(a))=\ld(a/K)$. Let
$b=(a,\si(a),\ldots,\si^n(a))$; then
$\ld(a/K)=\ld(b/K)=\mu(\si(b)/K(b))$. Furthermore, $b$ is a generic of
the difference subgroup $H_n$ (of $G\times \cdots\times G^{\si^n}$)
defined by 
$$H_n=\{(g_0,\ldots,g_n)\mid g_0\in H,
\bigwedge_{i=1}^ng_i=\si(g_{i-1})\}.$$
As $b=g(a)$ for some isomorphism $g:H\to H_n$ given by tuples of difference
polynomials, it suffices to prove the result for $b$ and $H_n$. 

Hence,
replacing $a$ by $b$, $H$ by $H_n$, we 
may assume that
$\si(a)\in K(a)^{alg}$ and $\ld(a/K)=\mu(\si(a)/K(a))$. Without loss of
generality, $H$ is Zariski dense in $G$, so that 
$G$ is connected.  

Let $\ell_0$ be defined as in Lemma \ref{sc-lem1}, take 
$\ell\geq \ell_0$, and consider the algebraic groups $U_\ell$,
$V_\ell$, where $U_\ell$ is the algebraic locus of
$(\si^{-\ell}(a),\si^\ell(a))$ over $K^{inv}$, and $V_\ell$ the algebraic
locus of $(a,\si^{-\ell}(a), \si^\ell(a))$ over $K^{inv}$. Then $V_\ell$ is an
algebraic subgroup of $G\times U_\ell$, and its images under the
projections $\pi_1:G\times U_\ell\to G$ and $\pi_2:G\times U_\ell\to U_\ell$
equal $G$ and $U_\ell$ respectively. 

We now apply 
Lemma \ref{lem5}, and use its notation and the notation  of Theorem
\ref{sc-thm1}. Note that $S_1$ is finite, so that in
particular $S_1$ 
is central in $G$ (since $G$ is connected). Let $f$ be the isogeny $G\to
G/S_1$, and $d=f(a)$. Then $d$ encodes the set $a+S_1$ of field conjugates of
$a$ over $K(\si^{-\ell}(a),\si^\ell(a))$. Recall that this field
contains the tuple $c$ of Theorem \ref{sc-thm1}, and that
$\mu(a/K^{inv}(c))=\mu(a/K^{inv}(\si^{-\ell}(a),\si^\ell(a)))$. By Proposition
\ref{prop2}(2), for some $r$ 
we obtain $\ld(f(a)^{p^r}/K^{inv})=\dd(a/K^{inv})$.

As constructed, our element $f(a)^{p^r}$ is in $K(a)_\si^{inv}$, not necessarily in
$K(a)_\si$. But for some $m$, $\si^m f(a)^{p^r}\in K(a)_\si$. 
We let $H'$ be the $\si$-closure of $\si^m\circ \Frob^r\circ f(H)$ in
$(G/S_1)^{\si^m\circ\Frob^r}$. Then $\si^m(d^{p^r})$ 
is a generic of $H'$, and $\si^m\circ \Frob^r\circ f$ defines a group
homomorphism 
$H\to H'$. (Here $\Frob$ denotes the Frobenius automorphism). 

If $K$ is perfect, then $\ld(d/K)= \dd(a/K)$, so we may take $H'$ to be
the $\si$-closure of $\si^m\circ  f(H)$.

\section{Comparison and/or equivalence of the various settings}

In this section we  first recall Willis' definitions and results on
totally disconnected locally compact groups (see \cite{W1}, \cite{W2}) and explain how they give
 our results for 
difference fields of  characteristic $0$. We  then compare the two
sets of results, in the group case and in the field case; and exhibit
some interesting translations. We end the section with the proof that any
totally disconnected locally compact  group is the inverse limit of
 automorphism groups of strongly minimal structures.

\para\vlabel{Wil1} {\bf The scale of a totally disconnected locally
  compact group}. 
Let $G$ be 
a totally disconnected locally compact
group, with 
a continuous automorphism $\alpha$. Let $U$ be an open compact subgroup
of $G$, and define
$$U_+=\bigcap_{n\in\nat}\alpha^n(U),\quad
U_-=\bigcap_{n\in\nat}\alpha^{-n}(U).$$
Say that $U$ is {\em tidy for $\alpha$} if it satisfies
\begin{itemize}
\item[{\bf T1}] $U=U_+U_-=U_-U_+$, and
\item[{\bf T2}] $\bigcup_{n\in\nat}\alpha^n(U_+)$ and
  $\bigcup_{n\in\nat}\alpha^{-n}(U_-)$ are closed in $G$. 
\end{itemize}
One then defines the {\em scale function of $\alpha$ on $G$} by 
$$s_G(\alpha)=[\alpha(U):\alpha(U)\cap U],$$
where $U$ is a tidy subgroup. That tidy subgroups exist and that the
scale function is well-defined is shown in \cite{W1}, Theorems 1 and 2. 

\smallskip
Let us now go to difference fields and see how the duality works. For
simplicity of notation we will assume that the characteristic is $0$; in
positive characteristic, analogous results are obtained if one replaces
everywhere the degree of a field extension by its separable degree. 
Let $K=\si(K)$ be a difference subfield of $\calu$, $a$ a tuple in
$\calu$ such that 
$\si(a)\in K(a)^{alg}$, and $L=K(a)^{alg}$. Set $$G=\aut(L/K), \qquad
V=\aut(L/K(a)).$$ Then $G$ is locally compact, and $V$ is a compact open
subgroup which is profinite. The action of $\si$ on $L$ induces a
continuous 
action $\alpha$ on $G$: 
$$\tau\mapsto \si\tau\si\inv,$$ which maps $V=\aut(L/K(a))$
onto $\aut(L/K(\si(a)))$. Then $V_+=\aut(L/K(a)_\si)$, and
$V_-=\aut(L/K(a)_{\si\inv})$ (where 
$K(a)_{\si\inv}=K(\si^{-n}(a), n\in\nat)$).

Condition T1 then corresponds to
$\mu(\si(a)/K(a))=\ld(a/K)$. Condition T2 is not so clear,
until one inspects Lemma~3(a) of
\cite{W1}: $\bigcup_{n\in\nat}\alpha^n(U_+)$ is closed if and only if 
$\bigcup_{n\in\nat}\alpha^n(U_+)\cap U=U_+$. This implies that
$\alpha^\ell(U_+)\cap 
U\subseteq U_+$ for $\ell>0$ and, assuming T1, a moment's thought shows
that it gives 
$\alpha(U)\supseteq U\cap 
\alpha^{\ell}(U)$. Thus, if $V$ is tidy, this tells us that $\si(a)\in
K(a,\si^\ell(a))$.

Thus, in characteristic $0$, the existence of tidy subgroups of $G$
together with this lemma give us (almost) Theorem~\ref{sc-thm1}. Indeed,
Theorem~1 of \cite{W1} gives a tidy subgroup $U$ which is compact open,
and 
therefore commensurable with $V$. I.e., if $K(b)$ is the subfield of $L$
fixed by $V$ then $K(a,b)$ is a finite extension of $K(a)$ and of
$K(b)$. However, inspection of the construction of this subgroup $U$
(see e.g. Lemma 3.3 in \cite{W2}) shows that it contains (a finite
intersection of transforms of) $V$. I.e.,
$b\in K(a)_\si$.  

The fact that an element which satisfies $\ld=\dd$ must also satisfy the
conclusion of Theorem~\ref{sc-thm1} is fairly clear, so the existence of
tidy subgroups led us to look closely at the proof of Theorem 1 of
\cite{W1} and to  discover the above mentioned implication of 
Lemma~3(a). It suggested that the result might be true in all
characteristic, 
but for that we needed to find a proof slightly more precise. 
We got more help  from Willis' definition of the group
$\call$ (see \cite{W1} page 347), which suggested that the field $F$ of
\ref{sc-thm1} might be large. However, the rest of our proof is somewhat
different from Willis'.

\para\vlabel{Wil2}{\bf Comparison of the results in the group and in the field context}. Below we will give a dictionary of how the various results relate to
each other. We first list the group-theoretic result (g), then immediately
below its field 
analogue (f). Many results are very similar, some are unexpected.\\[0.1in] 
(1)(g) The scale function does not depend on the chosen tidy
subgroup (Theorem~2 and/or Lemma~10 of \cite{W1}).\\
(f) Lemma \ref{lem4}
tells us that $\dd(a/K)$ is an invariant of the difference field
extension $K(a)_\si^{alg}/K$. See also
\ref{sc-lem1}(6): if $c,d$ satisfy $\ld=\dd$, then so does $(c,d)$. \\[0.2cm]
(2)(g) The modular function $\Delta(\alpha)$ of $\alpha$ equals
$s(\alpha)s(\alpha\inv)\inv$ (Corollary ~1 
 of \cite{W1})\\
(f) If $a$ and $b$ are equi-algebraic over $K$, then
${\ld(a/K)\over \ild(a/K)}={\ld(b/K)\over \ild(b/K)}$
(Lemma~\ref{lem4}(3)). \\[0.2cm] 
(3)(g) $s(\alpha^n)=s(\alpha)^n$ for $n>0$ (Corollary~3 of  \cite{W1}).\\
(f) $\dd_{\si^n}(a/K)=\dd(a/K)^n$ (Lemma \ref{sc-lem1}(6)).\\[0.2cm]
(4)(g) If $U$ is tidy for $\alpha$, and $\beta$ is conjugation by some element
$\tau\in U$, then $U$ is tidy for $\alpha \beta$, and $s(\alpha
\beta)=s(\alpha)$  (Theorem 3 of \cite{W1}, p.~356).\\
(f) This one is totally 
unexpected on the field side. Translated, it becomes: 
$$\hbox{If }
\ld(a/K)=\dd(a/K)\hbox{ and }\tau\in \aut(L/K(a)),\hbox{ then }
\ld_{\si\tau}(a/K)=\dd_{\si\tau}(a/K)=\dd(a/K).$$
This is a direct consequence of the following striking result, inspired
by the proof given in \cite{W1}: \\
{\bf Proposition}. If $a$ satisfies $\mu(\si(a)/K(a))=\ld(a/K)$, and $\tau\in \aut(K(a)^{alg}/K(a))$, then the difference fields
$(K(a)_\si,\si)$ and $(K(a)_{\si\tau},\si\tau)$ are isomorphic (by a
$K$-isomorphism taking $a$ to $a$). 

\prf Observe first that if $\rho_1,\rho_2\in\aut(K(a)^{alg}/K(a))$, then
the linear disjointness of 
$K(a)_{\si\inv}$ and 
$K(a)_\si$ over $K(a)$ implies the linear disjointness of
$\rho_1(K(a)_{\si\inv})$ and $\rho_2(K(a)_\si)$ over $K(a)$. In
particular, there is $\rho\in\aut(K(a)^{alg}/K(a))$ which 
agrees with $\rho_1$ on $K(a)_{\si\inv}$ and with $\rho_2$ on
$K(a)_\si$. 

One  shows by induction on $n$,
that  $K(a,\si(a),\ldots,\si^n(a))\simeq
K(a,\si\tau(a),\ldots,(\si\tau)^n(a))$ by a 
$K$-isomorphism (of fields) $f_n$ which sends $\si^i(a)$ to $(\si\tau)^i(a)$ for
$0\leq i\leq n$. For $n= 1$, $\tau\si\inv$ sends $(a,\si(a))$ to $(a,\tau(a))$.
 Assume given $f_n$, and observe that the field
$K((\si\tau)\inv(a),a)$ is precisely the image by $\tau\inv$ of the
field $K(\si\inv(a),a)$; we let $f_{-1}$ denote the restriction of
$\tau\inv$ to $K(\si\inv(a),a)$. By the remark above, and because
$f_{-1}$ and $f_n$ are the identity on $K(a)$, there is an element
$\rho\in \aut(K(a)^{alg}/K(a))$ which extends 
$f_{-1}\cup f_n$. Let $f_{n+1}$ be the restriction of $(\si\tau)\rho
\si\inv$ to $K(a,\ldots,\si^{n+1}(a))$.  \\[0.2cm]
%
%
(5)(g) $s(\alpha)=\min\{[\alpha(U):U\cap \alpha(U)]\mid U\hbox{ compact
  open}\}$; $[\alpha(U):\alpha(U)\cap U]=s(\alpha)\iff
[\alpha\inv(U):\alpha\inv(U)\cap U]=s(\alpha\inv)$  (Theorem~3.1 and
Corollary~3.11 of \cite{W2}).\\ 
(f) The equivalence of items (i) and
(ii) of Proposition~\ref{prop2}(1), and their equivalence with the
statement for $\si\inv$.  \\[0.2cm]
(6)(g) Let $H$ be a closed subgroup of $G$ such that $\alpha(H)=H$. Then
there is 
a tidy subgroup $U$ of $G$, such that $U\cap H$ is tidy for
$\alpha\rest_H$; furthermore $s(\alpha\rest_H)\leq s(\alpha)$ (Corollary
~4.2 and 
Proposition~4.3 of \cite{W2}).\\
(f) Let $M$ be a difference subfield of $L$
containing $K$. If $\ld(a/K)=\dd(a/K)$, then $\ld(a/M)=\dd(a/M)$: this
is clear using \ref{prop2}(1); $\dd(a/M)\leq \dd(a/K)$ is
obvious. However, Example 6.4 of \cite{W2} tells 
us that this is not the exact analogue of the group statement. \\[0.2cm]
(7)(g) Let $H$ be a closed normal subgroup of $G$ satisfying $\alpha(H)=H$, and
$\dot\alpha$ the automorphism of $G/H$ induced by $\alpha$. Then
$s(\alpha\rest_H)s(\dot\alpha)$ divides $s(\alpha)$ (Proposition 4.7 of
\cite{W2}).\\ 
(f) $\dd(a,b/K)\geq \dd(a/K(b)_\si)\dd(b/K)$ (Proposition
\ref{prop1}(1)). Thus we get a weaker result, but also under weaker
assumptions. On the other hand $\aut(L/K^{alg})$ has no proper closed
normal subgroup.\\[0.2cm]%
\para\vlabel{Wil3}{\bf Additional remark and results}.  We conclude with a remark on
some ingredients of our proof. We 
constantly use equation \ref{propld}($\#$), it is easy to derive the
analogue in the group context. The other ingredient we are using is the
tuple $c$ which encodes the tuple of minimal polynomials of $a$ over
a given field, see \ref{minpol}; its existence and properties guarantee
that certain infinite intersections are large. The analogue in the
group context 
exists, and can be stated as follows:\\

\hangafter=0\hangindent=20pt
\noindent Let $U$ be a compact open subgroup, $V$ a  compact
  subgroup of $G$, such that $[V:V\cap U]=N<\infty$. There is a compact open
   subgroup $W$ of $G$ which contains $V$, satisfies $[W:W\cap U]=N$, and
   contains all  subgroups with these properties. \\

\hangafter=0\hangindent=0pt
\noindent
This result is not difficult to prove, here is a sketch. Let $\calw$ be
the family of compact subgroups of $G$ which contain $V$ and satisfy
$[W:W\cap U]=N$. Note that this last condition is equivalent to $W\cdot
U=V\cdot U$ (where $W\cdot U$ denotes $\{wu\mid w\in W,u\in U\}$). The
family $\calw$ is non-empty ($V\in\calw$); observe that 
if $W_1,W_2\in\calw$, so does $W_1\cap W_2$, and therefore also $\langle
W_1 W_2\rangle$: this follows easily from $W_1\cdot W_2\cdot
U=W_1\cdot(W_1\cap W_2)\cdot 
U=W_1\cdot U$. Also, the closure of an element of $\calw$ is in $\calw$,
and this implies that $\calw$ has a unique maximal element, say
$W_0$. As  $\bigcap_{v\in V}v\inv U v$ is an open subgroup
which is  normalized
by $V$, it is contained in $W_0$, and therefore $W_0$ is open compact.

\smallskip 
When translated, our proof gives
a slightly different proof of the result in the group situation. Note
 the alternate definition of the scale function as
$$s(\alpha)=\lim_{k\rightarrow+\infty}[\alpha^k(U):U\cap\alpha^k(U)]^{1/k},$$
where $U$ is any compact open subgroup of $G$, and which comes from the
analogue of Lemma~\ref{lem4}(1). (This fact  was already
observed by R. G. M\"oller, \cite{Moe}.)
One can
also easily obtain 
the result corresponding to \ref{prop1}(7): \\

\hangafter=0\hangindent=20pt
\noindent If $U$ 
satisfies T1, and $W$ is a compact open subgroup which contains
$\alpha^{-\ell}(U)\cap 
\alpha^\ell(U)$ for all $\ell\gg 0$, then $W$ is tidy.\\

\hangafter=0\hangindent=0pt
\noindent
These results does
not seem to appear in either \cite{W1} or \cite{W2}. 

\para\vlabel{sm}{\bf Totally disconnected locally compact groups and strongly minimal sets}

If $T$ is a  disintegrated  strongly minimal theory\footnote{Recall that
  a theory $T$ is {\em strongly minimal} iff in any model $M$ of $T$, every
  definable subset of $M$ is finite or cofinite. It is
  {\em disintegrated} iff for any $A\subset M$, one has
  $\acl(A)=\bigcup_{a\in A}\acl(a)$.}, and $M$ is a model of
$T$,  then for any non-algebraic singleton $a\in M$, the group $\aut( \acl(a)
/ \acl(\emptyset))$ has the natural structure of  a totally disconnected locally compact
group (basic open sets are translates of stabilisers of finite sets;
note that $\aut(\acl(a)/a)$ is profinite and therefore compact).  Conversely, we will now
explain why any 
totally disconnected locally compact
group $G$ is a projective limit of ones that arise in this way.

Let $O$ be an open compact subgroup of $G$, and let $N_O$ be the
intersection of all conjugates of $O$. If $O'$ is an open subgroup of
$O$, then we have a natural onto map $G/N_{O'}\to G/N_O$, and the
intersection of all subgroups $N_O$, $O$ open compact, is $1$, so that 
$$G=\lim_{\leftarrow} G/N_O.$$
We will show that each $G/N_O$ is the automorphism group of a strongly
minimal disintegrated  set. Without loss of generality,
$N_O=1$, i.e., $O$ contains no proper normal subgroup of $G$.

Let $X=G/O$,  with $n$-ary
relations $R_a=Ga$ for any $a=(a_1,\ldots,a_n) \in X^n$ and $n\in\nat$,
i.e. $R_a$ is the $G$-orbit of $a$.  So $G$ acts on $M=(X,R_a)_{a}$
automorphically, transitively, and faithfully
because $O$ contains no proper normal subgroup. Let $\bar O$ be the
image of $O$ in $X$.  As $G$ acts transitively on $X$, to show that the
homomorphism $G 
\to \aut(M)$ is 
surjective, it suffices to show that 
   $O \to \aut(M/\bar O)$ is surjective.
 To show that $O \to \aut(M/\bar O)$ is surjective,
since $O$ is compact it suffices to see that the image is dense.
Indeed if  $h \in \aut(M/\bar O)$ and $h(a )= b$ for two $k$-tuples
$a,b$ of $X$, then $(b,\bar O)$ must be in the orbit of $(a,\bar O)$ since they
have the same (quantifier-free) type; so
$ga=b$ for some $g \in G$ with $g\bar O=\bar O$, i.e. $g \in O$.

Now $M$ is strongly minimal and disintegrated since the automorphism
group is transitive, and for any basic relation $R=R_a$, for some $m$,
$R(\bar O,x_1,...,x_m)$ holds for only finitely many elements $x_1,...,x_m$; see \cite{[Iv]} and the references therein.

Each element $g$ of $G$ defines an automorphism $\alpha$ of $M$ (via the
natural action of $G$ on $X$) and the
corresponding 
action on $G$ (viewed as $\aut(M)$) is conjugation by $g$. Thus the analogues of Theorems
\ref{sc-thm1} and \ref{sc-thm2} for strongly
minimal sets would give us Willis' Theorems 1 and 2  for inner automorphisms of $G$
(since quotienting by $N_O$ is irrelevant). On the other hand, if $G$ is
totally disconnected locally compact, so is $H=G \semdp \langle \si\rangle$ for
any automorphism $\si$ of $G$, so that only considering inner
automorphisms is not a restriction.

\end{document}